\documentclass{svproc}
\usepackage{amsmath, amssymb} 

\title{ECM Factorization with QRT Maps}
\titlerunning{ECM Factorization with QRT Maps}
\author{Andrew N.W. Hone\thanks{Work begun on leave at 
School of Mathematics \&  
Statistics, University of New South Wales, 
NSW 2052, Australia.}} 
\authorrunning{Andrew N.W. Hone} 
\tocauthor{Andrew N.W. Hone}

\institute{
School of Mathematics, 
Statistics \& Actuarial Science, 
University of
Kent,  Canterbury CT2 7NF, UK, 
\email{A.N.W.Hone@kent.ac.uk}
}

\newcommand{\beq}{\begin{equation}}  
\newcommand{\eeq}{\end{equation}}  
\newcommand{\bea}{\begin{eqnarray}}  
\newcommand{\eea}{\end{eqnarray}}  
\newcommand\la{{\lambda}}

\newcommand\al{{\alpha}}   
\newcommand\be{{\beta}} 
\newcommand\tal{\tilde{\alpha}}  
\newcommand\tbe{\tilde{\beta}}  
\newcommand\om{{\omega}}

\newcommand\rd{{\mathrm{d}}}  

\newcommand\sn{{\mathrm{sn}}}
\newcommand\cn{{\mathrm{cn}}}
\newcommand\dn{{\mathrm{dn}}}
\newcommand\cd{{\mathrm{cd}}}

\newcommand{\F}{{\mathbb F}}

\newcommand{\Q}{{\mathbb Q}}
\newcommand{\Z}{{\mathbb Z}}

\newcommand{\Pro}{{\mathbb P}}

\newcommand\Pc{{\mathcal{P}}}

\begin{document} 

\maketitle

\date{}

%\dedicatory{}

\begin{abstract} %({\bf Regular Research Paper}.) 
Quispel-Roberts-Thompson (QRT) maps are a 
family of birational maps of 
the plane which provide the 
simplest discrete analogue of an 
integrable Hamiltonian system, and are 
associated with elliptic fibrations in 
terms of biquadratic curves. Each 
generic orbit of a QRT 
map corresponds to a sequence 
of points on an elliptic 
curve. In this preliminary study, 
we explore versions of the 
elliptic curve method (ECM) for 
integer factorization based on performing 
scalar multiplication of a point 
on an elliptic curve by 
iterating three different QRT maps 
with particular initial data.  
Pseudorandom number generation and other 
possible applications are briefly discussed. 
\end{abstract}

\section{Introduction}

Elliptic curves are a fundamental tool in modern cryptography. 
The abelian group structure on an elliptic curve  makes it suitable for versions of Diffie-Hellman key exchange
and ElGamal key encryption, as well as providing techniques for primality testing and integer factorization, 
among many other applications relevant to network security  
\cite{cp,koblitz,stinson,yan}. In this paper we consider an approach to  
integer factorization using elliptic curves.  

The elliptic curve method (ECM) due to Lenstra \cite{lenstra} is one of the most effective methods known for finding 
medium-sized prime factors of large integers, in contrast to trial division, Pollard's rho method, or the 
$p-1$ method, which quickly find small 
factors, or sieve methods, which are capable of finding very large prime factors. 
For factoring an integer $N$, the basic idea 
of the ECM is to pick (at random) an elliptic curve $E$ and a point $\mathcal{P}\in E$, then compute the scalar 
multiple $s\Pc=\Pc+\cdots+\Pc$ ($s$ times) in the group law of the curve, using arithmetic 
in the ring $\Z/N\Z$, take a rational function $f$ on $E$ 
with a pole 
at  the point $\mathcal{O}$ corresponding to the identity in the group $E$, and evaluate  $f(s\Pc)$ 
for some 
$s$ chosen as the largest prime power less than some fixed 
bound $B_1$, 
or as 
the product of all such prime powers. 
For certain choices of $E$ and $\Pc$,  this computation  may lead to an attempt to divide by a 
non-unit in the ring, resulting in a factor of $N$ being found. 

To be more precise, 
traditionally 
one starts with a Weierstrass cubic defined over $\Q$, 
which can be taken with integer coeffients as 
%\beq\label{weier} 
%E: \, %\qquad 
$$y^2=x^3+Ax+B, \qquad 
A,B\in\Z,
$$ 
%\eeq 
%(where ), 
so that arithmetic $\bmod \,N$ corresponds to working with the pseudocurve (or group scheme) $E(\Z/N\Z)$ consisting  of 
all $(x,y)\in  (\Z/N\Z)^2$ that satisfy the cubic equation %(\ref{weier}), 
together with $\mathcal{O}$, the point at infinity; 
but when $N$ is composite 
the group addition $\Pc_1+\Pc_2$ is not defined for all pairs of points 
 $\Pc_1,\Pc_2\in E(\Z/N\Z)$. 
Typically $f$ is taken to be the coordinate function $x$, 
and the method is successful if  computing 
the scalar multiple $s\Pc$  leads to an $x$-coordinate with a denominator $D$
which is not a unit in 
$\Z/N\Z$, such that $\gcd(D,N)>1$ is a non-trivial factor of $N$. 
When this fortunate 
occurrence arises, it indicates that there is a prime factor 
$p|N$ for which 
$s\Pc = \mathcal{O}$ in the group law of the 
bona fide elliptic curve $E(\F_p)$, 
which is guaranteed if $s$ is a multiple of the order 
$\#E(\F_p)$. 

The original description of the ECM was based on computations 
with affine 
coordinates for a Weierstrass cubic; 
%The basic outline of the ECM described 
computing the scalar multiple $s\Pc$
%above 
is now known as ``stage 1'' of the ECM, and there is a further 
``stage 2'', due to Brent, involving computing multiples $\ell s\Pc$ for small primes $\ell$ less than some 
bound $B_2>B_1$, but here we only focus on %the first 
stage 1. 
Improvements 
in efficiency can be made %obtained 
by using various types of 
projective coordinates and/or Montgomery curves %coordinates/curves 
(see chapter 7 in \cite{cp}). 
However, all of these approaches share an 
inconvenient feature of 
the addition law for  $\Pc_1+\Pc_2$ 
on a Weierstrass cubic, namely that 
%%it must be 
%%split  into separate cases: %in particular,  
the formulae
 for $\Pc_2=\pm \Pc_1$ or $\Pc_2=\mathcal{O}$ are different 
 from the generic case.

An important new development was the proposal 
of Bernstein and Lange \cite{bl} to 
consider a different model for $E$, namely the 
% an elliptic 
Edwards 
curve \cite{edwards} 
\beq\label{edwards} 
E_d: \qquad 
x^2 +y^2=1+dx^2y^2  
\eeq 
($d$ is a parameter), 
for which the addition law 
\beq\label{edaddn} 
(x_1,y_1)+(x_2,y_2) =\left( 
\frac{x_1y_2+y_1x_2}{1 +d x_1x_2y_1y_2}, 
\frac{y_1y_2-x_1x_2}{1 -d x_1x_2y_1y_2}
\right) 
\eeq 
has the advantage that it is also 
valid for a generic pair of 
points 
$\Pc_1,\Pc_2\in E_d$, even when $\Pc_1=\Pc_2$, so it can be used for doubling 
(following \cite{bl}, we have used a rescaled curve compared with the original version in 
\cite{edwards}). 
The fact that the addition  law 
(\ref{edaddn}) on $E_d$ is unified  in this sense 
is implicit in the classical 
addition formula for the Jacobi sine function (see chapter XXII in \cite{ww}, 
or chapter 22 in \cite{nist}),  
for we have 
$$ 
\sn (z+w) = \frac{\sn (z) \cd (w) + \cd (z) \sn (w)}
{ 1+k^2 \sn (z) \sn (w) \cd (z) \cd (w)}, 
$$ 
$$ 
\cd (z+w) = \frac{\cd (z) \cd (w) - \sn (z) \sn (w)}
{ 1-k^2 \sn (z) \sn (w) \cd (z) \cd (w)}, 
$$ 
using Glaisher's notation for the quotient 
$\cd(z)= \cn(z)/\dn(z)=\sn (z+K)$,   
with the complete elliptic integral $K=K(k)$ being 
a quarter-period of the Jacobi sine, 
which yields (\ref{edaddn}) when we parametrize the 
points on $E_d$ by 
\beq\label{edparam}
(x,y)=\big(\sn(z),\cd(z)\big)=
\big(\sn(z),\sn(z+K)\big)
\eeq 
 and identify $d=k^2$.

It was shown in \cite{bl} that, compared 
with the Weierstrass representation and 
its variants, 
the Edwards 
addition law 
gives more efficient formulae for 
computing an addition step 
$(\Pc_1,\Pc_2)\mapsto \Pc_1+\Pc_2$ 
or a doubling step 
$\Pc_1\mapsto 2\Pc_1$, both of which are required to obtain the scalar 
multiple $s\Pc$ in subexponential time $O(\log s)$ via an addition chain. 
The implementation EECM-MPFQ introduced in 
\cite{bblp} gains even greater efficiency by using twisted 
Edwards curves,  
with %premultiplying 
an extra parameter $a$ in front of the term $x^2$ on 
the left-hand side of (\ref{edwards}), 
and further 
optimizing the ECM in other ways, including the use of 
projective coordinates in $\Pro^2$, 
extended Edwards coordinates in $\Pro^3$, and choosing curves with large 
torsion. 

In this paper we %are interested in exploring potential 
explore implementations of the 
ECM using other models of elliptic curves, %$E$, 
which 
arise %naturally 
in the context of QRT maps,  
%In order to do so, we 
%will making use of  
%QRT maps are 
an 18-parameter family of birational maps of the plane %, which were 
introduced 
by Quispel, Roberts and Thompson  \cite{qrt1} to unify %provide a 
%unifying framework for 
diverse examples of maps and functional relations 
%that had previously appeared 
appearing in dynamical systems, statistical mechanics 
and soliton theory. 
A QRT map is one of the simplest examples of a discrete integrable 
system, being a discrete avatar of a Hamiltonian system with one 
degree of freedom, with an invariant function (conserved quantity) 
and an invariant measure (symplectic form) \cite{duistermaat}.

Each orbit of a QRT map corresponds
to a sequence of points $\Pc_0 +n\Pc$ on a curve of genus one,
and in the special case %where the initial point 
$\Pc_0 =\mathcal{O}$  
the orbit consists of the scalar multiples $n\Pc$, being %and 
closely related to an elliptic divisibility sequence (EDS) %one of Morgan Ward's 
\cite{ward}. 
Thus  we can implement the ECM by iterating a QRT map with a special choice 
of initial data, and performing all the arithmetic in      
$\Z/N\Z$.

A terse overview of QRT maps is provided in the next section;  see \cite{duistermaat,iatrou1,iatrou2,tsuda} for further details. 
%%the interested 
%%reader can refer   to the extensive literature on the topic, such as \cite{duistermaat,iatrou1,iatrou2,tsuda}, 
%%or   Duistermaat's book \cite{duistermaat}. 
Section 3 briefly introduces Somos sequences and related 
EDS, showing how three particular examples of QRT maps arise in this context, namely the Somos-4 QRT map, 
the Somos-5 QRT map, and the Lyness map. Each of the subsequent sections 4-6 is devoted 
to %giving further details of 
one of these three types of QRT map, 
including the doubling map that sends any point 
$\Pc_1\mapsto2\Pc_1$, %together 
and a corresponding version of the 
ECM. In section 7 we analyse the complexity of scalar multiplication, concentrating %mainly 
on the 
Lyness case  in %both 
projective %and affine 
coordinates, and the final section contains some conclusions.

\section{A brief review of QRT maps}\label{qrtsec}
 
A QRT map can be constructed from %by starting with 
a biquadratic curve, of the general  form 
\beq\label{biq}
F(x,y) := \sum_{i,j=0}^2 a_{ij}x^iy^j=0.   
\eeq  
For generic coefficients $a_{ij}$, this is a smooth affine curve, and with the inclusion of  
additional points at infinity it  lifts  to a smooth curve  in $\Pro^1\times \Pro^1$, %as is seen 
by introducing homogeneous 
coordinates $\big( (X:W), (Y:Z)\big)$ and setting $x=X/W$, $y=Y/Z$ to 
obtain a homogeneous equation of bidegree $(2,2)$, that is   
$$
\hat{F}(X,W,Y,Z)=W^2Z^2 F(X/W,Y/Z)=0;
$$
this curve 
is a double cover of $\Pro^1$ with four branch points, so %it 
has genus one 
by Riemann-Hurwitz.  
A biquadratic curve admits two simple involutions, namely the horizontal/vertical switches 
given by  
$$
\iota_h: \, (x,y)\mapsto (x^\dagger,y), \qquad 
\iota_v: \, (x,y)\mapsto (x,y^\dagger), $$  where $x^\dagger$ is the conjugate root 
of (\ref{biq}), viewed as a quadratic in $x$, and similarly for $y^\dagger$; 
%where the equation (\ref{biq}), together with 
the Vieta formulae for the sum/product of the roots 
of a quadratic allow explicit birational expressions to be given for these two involutions. On a given 
biquadratic curve,  the QRT map is defined to be the composition of the two switches,  
%\beq\label{qrt} 
$$\varphi_{QRT}=\iota_v\circ\iota_h,$$ 
%\eeq 
%and general considerations about elliptic curves imply that 
which acts as a translation 
in the group law of the curve, $\varphi_{QRT}: \, \Pc_0\mapsto \Pc_0 +\Pc$, where 
the shift $\Pc$ is independent of the choice of initial point    $\Pc_0$ on the curve. 

So far the map $\varphi_{QRT}$ is restricted to a single curve, but to define a %QRT 
map 
on the plane one should  allow
each coefficient $a_{ij}=a_{ij}(\la)$  to be a linear function of a parameter $\la$, 
so that (\ref{biq}) becomes   
%a pencil of biquadratic curves, 
a biquadratic pencil,
\beq\label{pencil}
E_\la: \qquad F(x,y) \equiv F_1(x,y) +  \la\, F_2(x,y) =0. 
\eeq  
The map $(x,y)\mapsto \la =-F_1(x,y)/F_2(x,y)$, 
obtained by solving  (\ref{pencil}) for $\la$, 
defines an elliptic fibration of the plane over $\Pro^1$ (except at  finitely many base points 
where $F_1=F_2=0$). Each value of $\la$ corresponds to  
a unique curve in the pencil, where 
the map 
$\varphi_{QRT}$ is defined, and on each such curve, 
% and %by using 
a suitable 
combination of Vieta formulae %one obtains 
yields a birational expression which is 
independent of $\la$, so defines a birational map on the $(x,y)$ plane, 
also denoted $\varphi_{QRT}$. %Moreover, 
By construction the function $-F_1/F_2$ is constant on each orbit, so 
is a conserved quantity for the map $\varphi_{QRT}$ in the plane. 

Henceforth we restrict to the symmetric case $F(x,y)=F(y,x)$, so that 
each curve in the pencil also admits the involution   
$$
\iota:\,
(x,y)\mapsto (y,x),
$$    
%which makes 
making the horizontal/vertical switches conjugate to one another; 
%implies that 
thus $\varphi_{QRT}$ is a perfect square: 
$\iota_v=\iota\,\circ\,\iota_h\,\circ\,\iota$, hence $\varphi_{QRT}=(\iota\,\circ\,\iota_h)^2 =\varphi\,\circ\,\varphi$, 
where the ``square root'' of $\varphi_{QRT}$ is the symmetric QRT map 
%\beq\label{symqrt} 
$$\varphi = \iota\circ\iota_h.$$  
%\eeq 

As a simple example, note that the Edwards curve (\ref{edwards}) is a symmetric biquadratic, and 
we can identify $d=\la$ as the parameter of the pencil. Then %from 
the Vieta formula for the sum 
of the roots gives an expression that is independent of this parameter, 
and the symmetric QRT 
map $\varphi=\varphi_{Edwards}$ associated with this pencil has the very simple form 
$$ 
\varphi_{Edwards}: \qquad (x,y)\mapsto (y,-x),
$$
which is periodic with period four, i.e.\ $(\varphi_{Edwards})^4=\mathrm{id}$. 
This is another manifestation of the well known fact that Edwards curves have 4-torsion, 
or of the fact that the complete elliptic integral $K$ in   (\ref{edparam}) is a quarter-period 
of the Jacobi sine.

A generic symmetric QRT map is far from being so simple: starting from an initial point $\Pc_0$ in the plane, each orbit is 
a sequence of points $\Pc_n=\Pc_0 +n\Pc$ on a particular curve $E_\la$, 
and in general (at least, over an infinite field) the shift $\Pc$ need not be a torsion point.  
Even over a finite field $\F_p$, where every point is torsion, the order of 
$\Pc$ typically varies with the choice of curve in the pencil, 
i.e.\ with the value of $\la$.

In the cases of interest for the rest of the paper, the symmetric QRT map $\varphi$ can be written in 
multiplicative form, %which means 
so that the sequence of 
points  $\Pc_n$ %=\Pc_0 +n\Pc$  
has coordinates $(x,y)=(u_n,u_{n+1})$, where 
$u_n$ satisfies a recurrence of second order, %of the form 
\beq\label{qrtrec}
u_{n+2}\, u_n =\mathrm{R}(u_{n+1}),  
\eeq 
for a certain rational function $\mathrm{R}$ of degree at most two, with 
coefficients that are independent of $\la$ (cf.\ Proposition 2.5 in \cite{hones5}, 
or %see 
\cite{iatrou1,iatrou2} for
more details). 
%a detailed discussion of normal forms for QRT maps). 

%By considering three particular Somos sequences, and three particular 

\section{Somos and elliptic divisibility sequences} \label{seds}

A Somos-$k$ sequence satisfies a quadratic recurrence of the form 
\beq\label{somosk} 
\tau_{n+k}\tau_n = \sum_{j=1}^{\lfloor k/2\rfloor} \al_j \, \tau_{n+k-j}\tau_{n+j},  
\eeq 
where (to avoid %certain 
elementary cases)  it is assumed that $k\geq 4$ with at least two 
%non-zero 
parameters $\al_j\neq 0$. 
%(so that  two or more  terms appear on the right-hand side).  
It was a surprising empirical observation of Somos \cite{somos} that such rational recurrences
can sometimes generate integer sequences, which was proved by Malouf \cite{malouf} for 
the Somos-4 recurrence 
\beq\label{s4} 
\tau_{n+4}\tau_n = \al \, \tau_{n+3}\tau_{n+1}+\be \, (\tau_{n+2})^2,   
\eeq 
in the particular case that the coefficients are $\al=\be=1$ and the initial values 
are $\tau_0=\tau_1=\tau_2=\tau_3=1$ . A broader understanding 
came from the further observation that the recurrence (\ref{s4}) has the 
Laurent property \cite{gale}, %%: %in the sense that 
%%the iterates are Laurent polynomials 
%%in the initial values, defined over the ring of integer polynomials in the coefficients, 
that is 
$
\tau_n \in \Z [\al, \be, \tau_{0}^{\pm 1},  \tau_{1}^{\pm 1},\tau_{2}^{\pm 1},\tau_{3}^{\pm 1}]
$ %\quad 
$\forall n\in\Z$. Somos sequences arise from mutations 
in cluster algebras \cite{fordy_marsh} or LP algebras \cite{lp}, 
and as reductions of the bilinear discrete KP/BKP equations, %(the octahedron/cube recurrences), 
being  associated with translations on  Jacobian/Prym varieties for the spectral curve of a 
corresponding Lax matrix \cite{s6prym,hkw}.

The three simplest non-trivial examples of Somos recurrences, with two terms on the right-hand side, 
are the Somos-4 recurrence (\ref{s4}), the Somos-5 recurrence  
\beq\label{s5} 
\tau_{n+5}\tau_n = \tal \, \tau_{n+4}\tau_{n+1}+\tbe \, \tau_{n+3}\tau_{n+2},   
\eeq 
and the special Somos-7 recurrence 
\beq\label{s7} 
\tau_{n+7}\tau_n = a \, \tau_{n+6}\tau_{n+1}+b \, \tau_{n+4}\tau_{n+3}.    
\eeq 
%We have introduced additional coefficients $\tal,\tbe,a,b$ to distinguish between these different examples. 
All  three of them can all be reduced to two-dimensional maps of QRT type, hence 
their orbits correspond to sequences of points $\Pc_0+n\Pc$ on curves of genus one. 
(In contrast, generic Somos-6 sequences and %all other %types of 
Somos-7 sequences 
are associated with %sequences of 
points on Jacobians of genus 2 curves \cite{s6prym}.) 

To see the connection with  QRT maps, in (\ref{s4}) one should 
substitute %make the substitution 
%\beq\label{s4subs} 
\beq\label{s4urec}
u_n = \frac{\tau_{n-1}\tau_{n+1}}{\tau_n^2} \implies    u_{n+2}\,u_n=\frac{\alpha \, u_{n+1} + \beta}{(u_{n+1})^2},
\eeq 
%which yields the 
yielding 
a second order recurrence 
%The latter
which  can be reinterpreted as 
the map 
$$
(u_n,u_{n+1}) \mapsto (u_{n+1},u_{n+2})$$ 
in the plane, %%corresponding to 
and it turns out to be a symmetric QRT map; 
%of the form (\ref{symqrt}); 
for the associated biquadratic pencil and 
other details, see %are presented in 
section \ref{s4sec}. 
Similarly, for the Somos-5 recurrence (\ref{s5}), one can make  the substitution   
\beq\label{s5urec} 
u_n = \frac{\tau_{n-2}\tau_{n+1}}{\tau_{n-1}\tau_{n}}\implies    
u_{n+2}\,u_n=\frac{\tal \, u_{n+1} + \tbe}{u_{n+1}},  
\eeq  
%substituted in (\ref{s5}) produces the recurrence 
%\beq\label{s5urec}u_{n+2}\,u_n=\frac{\tal \, u_{n+1} + \tbe}{u_{n+1}},  %\eeq 
where the latter recurrence for $u_n$ 
is equivalent to the QRT map described in section \ref{s5sec}. 
Finally, for the special Somos-7 recurrence (\ref{s7}) one should  substitute %take 
\beq\label{lurec} 
u_n = \frac{\tau_{n-3}\tau_{n+2}}{\tau_{n-1}\tau_{n}}\implies  
u_{n+2}\,u_n=a \, u_{n+1} + b,   
\eeq  
reducing the order from seven to two.  
%\beq\label{lurec}u_{n+2}\,u_n=a \, u_{n+1} + b, \eeq 
The %latter 
recurrence for $u_n$ in (\ref{lurec}) 
is known in the literature as the Lyness map, after the 
particular periodic case $b=a^2$ found in \cite{lyness}; 
for details see section \ref{lsec}.
%Although the substitutions in (\ref{s4subs}), (\ref{s5urec}) and (\ref{lurec}) 
The first two of these substitutions were derived in an ad hoc way in \cite{honeblms} and \cite{hones5},  
%for the first two, 
but they all have a very natural interpretation in the theory of cluster 
algebras \cite{fordyhone}, which implies %allowing one to deduce 
that these are the only Somos-$k$ 
recurrences that  can be reduced to two-dimensional maps.

Morgan Ward's elliptic divisibility sequences (EDS) \cite{ward} are 
sequences of integers $\tau_n$ with  
$\tau_0=0$, $\tau_1=1$, $\tau_2,\tau_3, \tau_4\in\Z$ and $\tau_2|\tau_4$,
subject to the relations 
\beq\label{eds1} 
\tau_{n+m}\tau_{n-m} = (\tau_{m})^2 \tau_{n+1}\tau_{n-1} -\tau_{m+1}\tau_{m-1} (\tau_{n})^2, 
\eeq  
\beq\label{eds2} 
\tau_2\tau_{n+m+1}\tau_{n-m} = \tau_{m+1}\tau_m \tau_{n+2}\tau_{n-1} -\tau_{m+2}\tau_{m-1} \tau_{n+1}\tau_n
\eeq  
for all $m,n\in\Z$. 
An EDS corresponds to a sequence of points   $n\Pc$ on an elliptic curve 
%defined 
over $\Q$. 
%If $\tau_2,\tau_3, \tau_4$ are all non-zero, then together with $\tau_1$ they determine the rest of the sequence, since 
The relation (\ref{eds1}) for $m=2$ is a special case of the Somos-4  %sequence, 
recurrence (\ref{s4}), with  
%with coefficients specified by two of the initial values as  
$\al=(\tau_2)^2$, $\be = -\tau_3$; similarly 
(\ref{eds2}) with $m=2$ gives a special case of (\ref{s5}), 
and a linear  combination of this relation for $m=2$ and $m=3$ yields 
(\ref{s7}) with the coefficients/initial values related in a  
particular way. 
The fact that the same EDS satisfies these higher Somos relations \cite{swartvdp} 
%connections between these special sequences 
provides one way to derive the isomorphisms between the associated biquadratic curves and 
a Weierstrass cubic  in Theorem \ref{iso} below, which can also be deduced  from results in \cite{swahon}.

\section{Somos-4 QRT map} \label{s4sec} 

Here %and in the subsequent two sections 
we give further details of the QRT map 
defined by (\ref{s4urec}) and the associated family of curves. 
%\noindent 
\begin{align}\label{s4qrt}
&\mathbf{QRT\,\, map:} \qquad \qquad 
\varphi: \, 
(x,y)  \mapsto 
\left(  y , \big({\alpha \, y + \beta}\big)/({xy^2})\right)  \\ 
%%\beq\label{s4qrtinv}%%\begin{array}{l}%%\mathbf{Inverse \,\,map:} \qquad 
%%\varphi^{-1}: \quad (x,y)  \mapsto \left(  \frac{\alpha \, x + \beta}{x^2y},x\right) . \end{array}\eeq 
\label{s4curve} 
&\mathbf{Pencil\,\, of\,\, curves:}\qquad 
x^2y^2+\al \, (x+y) + \beta -J\, xy=0. \\% \qquad\qquad\qquad
\label{einv} 
&\mathbf{Elliptic \,\,involution:} \quad 
\iota_E:\, 
(x,y)\mapsto \big(x,(\al \,x+\beta)/(x^2y)\big). \\ %\qquad\qquad
\label{idsh} %\begin{array}{l}
&\mathbf{Identity\,\,element \,\, \& \,\, shift:} %\,\, multiples:}\qquad \qquad
\qquad  
\mathcal{O}=(\infty,0), \quad \mathcal{P}=(0,-\beta/\alpha).  \\  
\nonumber 
&\mathbf{Doubling \,\,map:} \qquad\psi: \, (x,y)   \mapsto \\ 
\label{s4dub}
&\left( 
\frac{\al\,(x-y)y\,(\al \,x + \beta-x^3y)}{(\al \, x+\beta - x^2y^2)^2} , 
-\frac{(\al \, x+\beta - x^2y^2)(\al \, y+\beta - x^2y^2)}{\al \, xy(x-y)^2} 
\right) .  
\end{align}
The map  (\ref{s4qrt})  preserves the symplectic form 
$\om = (xy)^{-1} \rd x \wedge \rd y$, 
%\beq\label{omega} \om = \frac{1}{xy}\, \rd x \wedge \rd y, \eeq 
that is, 
$\varphi^*(\omega) = \omega$, 
and %it can be checked directly that 
the doubling map $\psi$ 
%also doubles the symplectic form, that is   
gives $\psi^*(\omega) = 2\omega$; the same is true for 
the Somos-5/Lyness maps.
Each orbit of  $\varphi$ 
lies on  a fixed biquadratic curve of the form (\ref{s4curve}), 
%corresponding to a particular value of $J$, 
with $\la=-J$ being the parameter of the pencil %in the form 
(\ref{pencil});  
equivalently, solving (\ref{s4curve}) for 
$J=J(x,y)$ gives 
%a function of $x,y$ which is 
a conserved quantity for the map.  
On any curve (\ref{s4curve}), the  elliptic involution (\ref{einv}) %which 
sends any point $\Pc \mapsto -\Pc$. % is %given by  $\iota_E:\, (x,y)\mapsto \big(x,(\al \,x+\beta)/(x^2y)\big)$. 
%We can construct 
A special sequence of points $(u_n,u_{n+1})$ on the curve 
is generated by iterating (\ref{s4qrt}) with a suitable starting point, 
corresponding to the scalar multiples $n\Pc$ of a particular point $\Pc$ (the shift).  
%In order to have an 
To have  both coordinates  finite 
and non-zero, one should start with  
%initial value for iteration of the affine map (\ref{s4qrt}), with 
\beq\label{s4init} 
2\mathcal{P} = ( -\beta/\alpha, -\al (\al^2+\beta J)/\beta^2)=(u_2,u_3).
\eeq  
However, in order to calculate a particular scalar multiple $s\Pc$ in time 
$O(\log s)$, rather than $O(s)$,  one must employ the doubling map on the curve, using some variant 
of the ``double-and-add'' method (an addition chain).

We can now present a version of the ECM based on the QRT map (\ref{s4qrt}). 

\vspace{.1in}
\noindent {\bf Algorithm 1: ECM with Somos-4 QRT.} 
\textit{To factorize %a composite integer 
$N$,  
pick %parameters 
$\alpha,\beta,J\in\Z/N\Z$ at random,   and some integer $s>2$.
Then starting from the point $2\Pc=(u_2,u_3)$ 
on the curve (\ref{s4curve}), given by (\ref{s4init}), use the QRT map  (\ref{s4qrt}) 
to perform addition steps, and
(\ref{s4dub}) to perform doubling steps, 
working in $\Z/ N\Z$, to compute $s\Pc=(u_s,u_{s+1})$. Stop if, for 
some denominator $D$,  $g= \gcd(D,N)>1$ 
appears at any stage; 
when $g<N$  the algorithm has been successful, but  if $g=N$ or no forbidden divisions appear then 
restart with new $\alpha,\beta,J$, and/or 
a larger value of $s$.} 

%%Strictly speaking, one also needs to check that 
%%$\gcd(\al,N)=1=\gcd(\be,N)$ at the start, otherwise 
%%$2\Pc$ is (\ref{s4init}) undefined, and the discriminant of (\ref{weier}),  
%%for $A,B$ as  in Theorem \ref{iso} below, 
%%requires $4A^3+27B^2$  coprime to $N$; if one of these conditions fails then 
%%a factor of $N$ has already been found. % without implementing the algorithm. 
%%Also, rather than picking $J$ it is simpler to pick  the $y$-coordinate of $2\Pc$ in (\ref{s4init}) at random. 

\begin{example} \label{s4exa} 
Given $N=1950153409$ we pick $\al=\be=1$, $J=4$ to find $(u_2,u_3)=(-1,-5)$, take $s=12$, and compute the sequence  $(u_n \bmod N)$, that is % which begins 
$$ \begin{array}{l}
\infty, 0, 
-1,
-5, 
1482116591, 
121884579, 
452175879,
1062558798, 
154165861, 
\\ 
1566968710,
1329544730, 
56956778,
\end{array} 
$$
where the last term is $u_{11}$; then $g=\gcd(u_{11},N)=16433$ so the algorithm terminates. 
Of course, not all the above terms are necessary, since by writing $12=2^2\times(2+1)$ it is 
more efficient to %enough to 
compute the addition chain $2\Pc\mapsto 3\Pc\mapsto 6\Pc\mapsto 12\Pc$ using (\ref{s4qrt}) and (\ref{s4dub})
as 
$$ 
(u_2,u_3)\overset{\varphi}{\mapsto} 
(u_3,u_4)\overset{\psi}{\mapsto}
(u_6,u_7)\overset{\psi}{\mapsto} ???
$$
and then observe that the denominator 
$\al x +\be -x^2y^2$ in (\ref{s4dub}) has common factor $g>1$ with $N$ when $(x,y)=(u_6,u_7)$. 
\end{example}

\section{Somos-5 QRT map}\label{s5sec}

Here we describe the features of the QRT map corresponding to recurrence (\ref{s5urec}).  
%is equivalent to the 
\begin{align} 
\label{s5qrt}
&\mathbf{QRT\,\, map:} \qquad \qquad
\varphi: \,
(x,y)  \mapsto 
\left(  y , \big({\tal \, y + \tbe}\big)/({xy})\right). \\ 
\label{s5curve} 
&\mathbf{Pencil\,\, of\,\, curves:}\qquad
xy(x+y)+\tal \, (x+y) + \tbe -\tilde{J}\, xy=0. \\ 
\nonumber 
&\mathbf{Elliptic \,\,involution:} \quad\iota_E:\quad (x,y)\mapsto (y,x). \\ 
\nonumber 
&\mathbf{Identity\,\,element \,\, \& \,\, shift:} \qquad \mathcal{O}=(\infty,\infty), \qquad   
\mathcal{P}=(\infty, 0). \\ 
\nonumber
&\mathbf{Doubling \,\,map:} \qquad\psi: \, (x,y)   \mapsto 
\end{align} 
%\\
\beq 
\label{s5dub}
%&
\left( 
\frac{(x^2y-\tal x-\tbe)(x^2y-\tal y-\tbe)}{x(x-y)(xy^2-\tal x-\tbe)}, 
\frac{(xy^2-\tal x-\tbe)(xy^2-\tal y-\tbe)}{y(y-x)(x^2y-\tal y-\tbe)}
\right) .
%\end{align}
\eeq  
%The invariant symplectic form $\om=(xy)^{- 1}\rd x\wedge \rd y$ is the same as  for the Somos-4 QRT map. 
The double  of the 
translation point (shift)  is  
%$\mathcal{P}=(\infty, 0)=(u_1,u_2)$, %\quad  
$2\mathcal{P} =(0,-\tbe/\tal)=(u_2,u_3)$, %$$ 
so %iteration of (\ref{s5qrt}) 
to obtain the sequence %$(u_n,u_{n+1})$ corresponding to the 
of multiples $n\Pc$ 
%by iterating %iteration of the affine map , 
one must start 
%should begin 
with %the finite 
\beq\label{s5init} 
3\mathcal{P} = ( -\tbe/\tal, \tilde{J}+\tal^2 /\tbe+\tbe/\tal)=(u_3,u_4).
\eeq

We can paraphrase Algorithm 1 to get another version of the ECM. 

\vspace{.1in}
\noindent {\bf Algorithm 2: ECM with Somos-5 QRT.} 
\textit{To factorize %a composite integer 
$N$,  
%modify the  approach used in Algorithm 1, by picking positive integer 
pick %parameters 
$\tal,\tbe,\tilde{J}\in\Z/N\Z$ at random,   and some integer $s>3$. 
Then starting from  $3\Pc=(u_3,u_4)$ 
on the curve (\ref{s5curve}),  given by (\ref{s5init}), use  (\ref{s5qrt}) 
to perform addition steps, and %either 
(\ref{s5dub}) %or (\ref{s5dubr}) 
to perform doubling steps, 
working in $\Z/ N\Z$, to compute $s\Pc=(u_s,u_{s+1})$. Stop if, for 
some denominator $D$,  $g= \gcd(D,N)$ 
with $1<g<N$ appears at any stage.} 

\section{Lyness map} \label{lsec}

The real and complex dynamics of the recurrence (\ref{lurec}), known as the Lyness map, has been studied by many authors,  %For 
with a very detailed account %of its geometrical properties, see 
in \cite{duistermaat}.   
\begin{align}
&\mathbf{QRT\,\, map:} \qquad
\varphi: \,\,  
(x,y)  \mapsto 
\left(  y , \frac{a \, y + b}{x}\right) . \label{lqrt}  \\
%\eeq 
&\mathbf{Pencil\,\, of\,\, curves:}\nonumber \\ 
\label{lcurve} 
& \qquad xy(x+y)+a\, (x+y)^2+(a^2+b) \, (x+y) + ab -K\, xy=0. 
\\ 
&\mathbf{Elliptic \,\,involution:} \qquad   \iota_E:\quad (x,y)\mapsto (y,x). \nonumber  
\\
%$$ 
%\beq
\label{lshift} 
&\mathbf{Identity\,\,element \,\, \& \,\, shift:}  \quad \mathcal{O}=(\infty,\infty), \quad 
\mathcal{P}=(\infty, -a). 
\end{align} 
%\eeq
%$$ 
\begin{align}
%\eeq 
%$$\begin{array}{r}
&\mathbf{Doubling \,\,map:} \quad\psi: \,\, (x,y)   \mapsto \Big(R(x,y), R(y,x)\Big), \nonumber \\ %\qquad\qquad\qquad\qquad\qquad\qquad\qquad\qquad\qquad\qquad
%\end{array}
%$$ 
%\beq
\label{ldub}
&R(x,y) = \frac{(xy-ay-b)(x^2y-a^2x-by-ab)}{x(x-y)(y^2-ax-b)}.
\end{align} 
%Once again, 
%The invariant symplectic form for the  QRT map is the same as before.  
%preserves the $\om=(xy)^{- 1}\rd x\wedge \rd y$, 
%and $\psi^*(\om)=2\om$.
%Solving (\ref{lcurve}) for $K$ gives a function of $x,y$ which is a conserved quantity, 
%and this gives a pencil of curves with parameter $\la=-K$. % for the map.   
Doubling and tripling $\Pc$ gives 
$2\mathcal{P} =(-a,0)$, %=(u_2,u_3)$, %\quad 
$3\mathcal{P} = ( 0, -b/a)$, %=(u_3,u_4)$
so to obtain the multiples $n\Pc=(u_n,u_{n+1}) $ by iteration of (\ref{lqrt}) and (\ref{ldub}) 
one should begin with 
\beq\label{linit} 
4\mathcal{P} = \left( -\frac{b}{a}, -a-\frac{b(Ka+b)}{a(a^2-b)}\right)=(u_4,u_5) 
.
\eeq 

Henceforth it will be assumed that $b\neq a^2$, since otherwise all orbits of (\ref{lurec}) are periodic with period five, 
meaning that $\Pc$ is a 5-torsion point on every curve in the pencil. This special case 
is the famous Lyness 5-cycle \cite{lyness}, 
related to the associahedron $K_4$ via the $A_2$ cluster algebra, and to the Abel pentagon identity for the 
dilogarithm \cite{nakanishi}, among many other things. 

The above formulae (and those for Somos-4/5) can all be obtained via the birational equivalence of 
curves  described in the following theorem (cf.\ \cite{swahon}). 

\begin{theorem} \label{iso} 
Given a fixed choice of rational point $\Pc=(\nu,\xi)\in\Q^2$ on a Weierstrass cubic 
$$E(\Q): \, 
(y')^2=(x')^3+Ax'+B$$  
over $\Q$,  a point $(x,y)$ on a Lyness curve (\ref{lcurve}) is given in terms of  
$(x',y')\in E(\Q)$ by %the formulae 
%\beq \label{xyxyprime}
$$x=- \frac{\be(\al u+\be)}{uv} -a, \quad   
y=-\be uv-a,$$  
%\eeq
where $$(u,v)=\left(\nu - x', \frac{4\xi y' + Ju-\al}{2u^2}\right) $$ are the coordinates of a point on the Somos-4 curve (\ref{s4curve}), 
and the parameters are related  by %the relations  
\beq\label{abK}
a=-\al^2 -\be J,\quad  b = 2a^2+a\be J-\be^3, \quad  K=-2a-\be J, 
\eeq 
with 
$$ 
\al=4\xi^2, \quad J=6\nu^2+2A, \quad \be = \frac{1}{4}J^2-12\nu\xi^2.
$$ 
Also, $$\left(-\frac{x+a}{\be},-\frac{y+a}{\be}\right)$$ 
is a point on the Somos-5 curve (\ref{s5curve}) with parameters 
$$ \tal =-\be, \quad \tbe = \al^2+\be J, \quad \tilde{J}=J.$$ 
Conversely, given $a,b,K\in\Q$, a point $(x,y)$ on  (\ref{lcurve}) corresponds  to $(\bar{x},\bar{y})\in\bar{E}(\Q)$, 
a twist of $E(\Q)$ with coefficients $\bar{A}=\al^2\be^4A$, $\bar {B}=\al^3\be^6B$, and $\Pc$ in (\ref{lshift}) 
corresponds to the %rational 
point $(\bar{\nu},\bar{\xi})=\big(\frac{1}{12}(\be J)^2-\frac{1}{3}\be^3,\frac{1}{2}\al^2\be^3)\in\bar{E}(\Q)$. 
\end{theorem}

\vspace{.1in}
\noindent {\bf Algorithm 3: ECM with Lyness.} 
\textit{To factorize %a composite integer 
$N$,  
%modify the  approach used in Algorithm 1, by picking positive integer parameters 
pick $a,b,K\in\Z/N\Z$ at random,   and some integer $s>4$. 
Then starting from  $4\Pc=(u_4,u_5)$ on the curve (\ref{lcurve}), given by (\ref{linit}), 
 use  (\ref{lqrt}) 
to perform addition steps, and %either 
(\ref{ldub}) %, (\ref{ldubr}) or (\ref{ldubr2})
 to perform doubling steps, 
working in $\Z/ N\Z$, to compute $s\Pc=(u_s,u_{s+1})$. Stop if, 
for 
some denominator $D$,   $g= \gcd(D,N)$ 
with $1<g<N$ appears  at any stage.}

%Comparison formulae: 

\section{Complexity of scalar multiplication} 
% in the Lyness case}

Of the three  symmetric QRT maps %considered 
above, the Lyness map (\ref{lqrt}) 
is the simplest, so we focus on that for our analysis. 
Before 
proceeding, we can make the simplification %, which is that we can set the 
%parameter 
$a\to 1$ without loss of generality, since over $\Q$ we can always rescale %the affine coordinates 
$(x,y)\to (ax,ay)$ and redefine $b$ and $K$.
To have an efficient version of Algorithm 3, it is necessary to work in projective 
coordinates, to avoid costly modular inversions; then only a single gcd 
need be calculated at the end. For cubic curves it is most common to work in 
the projective plane $\Pro^2$ (or sometimes Jacobian coordinates in the weighted projective space 
$\Pro (1,2,3)$ are used for Weierstrass cubics  \cite{gjm}).
However, for the biquadratic cubics (\ref{lcurve}), % it turns out that 
$\Pro^1\times \Pro^1$ is better, 
since %the expression for 
doubling with %map 
(\ref{ldub}) is of higher degree in $\Pro^2$. 

In terms of  projective coordinates in $\Pro^1\times \Pro^1$, the Lyness map (\ref{lqrt}) becomes 
\beq\label{lqrtp}  
\Big( (X:W), (Y:Z)\Big)\mapsto \Big( (Y:Z), ((aY+bZ)W:XZ)\Big) 
\eeq 
%where 
%(we  included 
%the parameter 
%$a$ for completeness). 
%Upon setting 
Then taking $a\to 1$,  
%for convenience then 
each addition step %, adding the point $\Pc$ 
using (\ref{lqrtp}) requires 
$2{\bf M}+1{\bf B}$, i.e.\ two multiplications and one multiplication by %the 
parameter $b$. 
 
The affine doubling map  (\ref{ldub})  %$\psi$ 
for the Lyness case %, given by the affine map (\ref{ldub}) with (\ref{lrfn}),  
lifts to the projective version  
\beq\label{ldubproj} 
\Big( (X:W), (Y:Z)\Big)\mapsto 
\Big( (A_1B_1:C_1D_1), (A_2B_2:C_2D_2)\Big), 
%\Big( (X^*:W^*), (Y^*:Z^*)\Big), 
\eeq 
where 
$$ 
X^*=A_1 B_1, \quad W^*=C_1D_1, \quad Y^*=A_2B_2, \quad Z^*=C_2D_2, 
$$ 
$$ 
A_1=A_++A_-, \quad A_2 =A_+-A_-, \quad B_1=B_++B_-, \quad B_2 =B_+-B_-,   
$$ 
$$
C_1=2XT, \quad C_2=-2YT, \quad D_1=ZA_2+C_2, \quad D_2=WA_1+C_1, 
$$ 
with $A_-=aT$ and 
$$ 
A_+=2G-aS-2H',  %\, A_-=aT, 
\,  B_+=S(G-a^2 H-H')-2aHH', 
\, B_-=T(G-a^2 H+H'), 
$$ 
$$ 
S=E+F, \,\,\, T=E-F, \,\,\, E=XZ, \,\,\, F=YW, \,\,\, G=XY, \,\,\, H=WZ, \,\,\, H'  =bH. 
$$ 
Setting $a\to 1$ once again for convenience, %as before 
and using the above formulae,  we see that doubling can be achieved 
with $15{\bf M}+1{\bf B}$. (To multiply by $2$ use  %can be  
%replaced by 
addition: $2X=X+X$.)  

%It is worth comparing 
This should be compared with %%the state of the art, 
 EECM-MPFQ \cite{bblp}: using twisted Edwards curves 
%%(\ref{edwards}), or rather their twisted versions 
$ax^2+y^2=1+dx^2y^2$ 
%%, to gain efficiency over %%compared with the standard the implementation of ECM with Weierstrass cubics. 
%% For  twisted Edwards curves in 
in $\Pro^2$
the projective  addition formula   requires 
$10{\bf M}+1{\bf S}+1{\bf A}+1{\bf D}$ (${\bf S},{\bf A}, {\bf D}$ denote squaring and
multiplication by the parameters $a,d$, respectively), while doubling only takes 
$3{\bf M}+4{\bf S}+1{\bf A}$. 
So the Lyness addition step (\ref{lqrtp}) %with Lyness curves 
is much more efficient than for twisted Edwards, but doubling 
requires twice as many multiplications. 
For any %what 
addition chain, %%is used 
the number of doublings will be $O(\log s)$, so employing Algorithm 3 to 
carry out the ECM with  the Lyness map 
in projective coordinates %, as above, 
should require on average roughly twice as many multiplications per bit as for    EECM-MPFQ.

\section{Conclusions} 

Due to the complexity of doubling, it  appears that scalar multiplication with Lyness curves is 
not competitive with the state of the art using twisted Edwards curves. However, in a follow-up study \cite{lynesseff},
we have shown that the projective doubling map (\ref{ldubproj}) for Lyness curves can be made efficient 
by distributing it over four processors in parallel, dropping the effective cost to $4{\bf M}+1{\bf B}$. 
On the other hand,  %but again 
this is still roughly twice the cost of the best known algorithm for doubling 
wih four processors on twisted Edwards curves in the special case $a=-1$ \cite{hisil}. 

However, by Theorem 
\ref{iso}, any elliptic curve over $\Q$ is isomorphic to a Lyness curve, while twisted Edwards curves only correspond to a subset of such curves. Thus there may be other circumstances, whether for the ECM or for 
alternative cryptographic applications, where 
Lyness curves and QRT maps will prove to be useful. For instance,  
one could use families of Lyness curves 
with torsion subgroups that are impossible with twisted Edwards curves in    EECM-MPFQ. 
Also, Bitcoin uses the 
curve $y^2=x^3+7$, known as secp256k1, which cannot be expressed in twisted Edwards form. 
    
The remarkable simplicity of the addition step (\ref{lqrtp}) means that it might also be suitable for pseudorandom 
number generation. In that context, it would be worth exploring non-autonomous versions 
of QRT maps $\bmod \,N$. 
For example, the recurrence 
\beq
\label{qp}u_{n+2}u_n =  u_{n+1}+b_nq^n, \qquad b_{n+6}=b_n 
\eeq 
is a $q$-difference Painlev\'e  version of the Lyness map (\ref{lurec}) (see \cite{hi}), 
%each step of which can be regarded as addition of a point $\Pc$ on an elliptic curve which changes with $n$.  
and over $\Q$ the arithmetic behaviour 
of such equations appears to be analogous to the autonomous case \cite{halburd}, 
with polynomial growth of logarithmic heights; 
although for (\ref{qp}) the growth 
is cubic rather than quadratic as in the elliptic curve case. It is interesting to compare this with the 
case where $q=1$ and the   coefficient $b_n$ is periodic with a period that does not divide 6, 
when generically (\ref{qp})  appears to display chaotic dynamics \cite{cima}, 
e.g.\  the period 5 example  $u_{n+2}u_n =  u_{n+1}+b_n$, $b_{n+5}=b_n$, 
for which the logarithmic height along orbits in $\Q$ grows exponentially with $n$.   
Working $\bmod \,N$, it would be worth carrying out a comparative study of  
the pseudorandom 
sequences generated by (\ref{qp}) to see how the behaviour for $q\neq 1$ differs  from 
the Lyness case (\ref{lurec}), and the effect of changing the period of $b_n$.

\section*{Acknowledgments} 
This research was supported by 
Fellowship EP/M004333/1  from the 
Engineering \& Physical Sciences Research Council, UK. The author thanks the 
School of Mathematics and Statistics, University of New South Wales, for 
hosting him twice during 2017-2019 as a Visiting Professorial Fellow, 
with funding from the Distinguished Researcher Visitor Scheme. 
He is also grateful to   
John Roberts and Wolfgang Schief for providing additional support during his 
stay in Sydney, where the idea behind this work originated, and to Reinout Quispel  
for useful discussions and hospitality 
during his visit to  Melbourne in May 2019.

\end{document}